\documentclass[11pt]{amsart}
\usepackage{mathrsfs}
\usepackage{graphicx}
\newtheorem{thm}{Theorem}[section]
\newtheorem{cor}[thm]{Corollary}

\theoremstyle{definition}

\theoremstyle{remark}

\numberwithin{equation}{section}

\begin{document}
\hfill{\emph{Annals of Operations Research}, 141 (2006), 19-52.}\\
\title[Multiserver retrial queues]{Analysis of Multiserver Retrial Queueing System: A Martingale Approach and an Algorithm of Solution}%
\author[Abramov]{Vyacheslav M. Abramov}%
\address{School of Mathematical Sciences, Monash University, Clayton, Victoria 3800, Australia}%
\email{Vyacheslav.Abramov@sci.monash.edu.au}%


\begin{abstract}
The paper studies a multiserver retrial queueing system with $m$
servers. Arrival process is a point process with strictly
stationary and ergodic increments. A customer arriving to the
system occupies one of the free servers. If upon arrival all
servers are busy, then the customer goes to the secondary queue,
orbit, and after some random time retries more and more to occupy
a server. A service time of each customer is exponentially
distributed random variable with parameter $\mu_1$. A time between
retrials is exponentially distributed with parameter $\mu_2$ for
each customer. Using a martingale approach the paper provides an
analysis of this system. The paper establishes the stability
condition and studies a behavior of the limiting queue-length
distributions as $\mu_2$ increases to infinity. As
$\mu_2\to\infty$, the paper also proves the convergence of
appropriate queue-length distributions to those of the associated
`usual' multiserver queueing system without retrials. An algorithm
for numerical solution of the equations, associated with the
limiting queue-length distribution of
retrial systems, is provided.\\



\end{abstract}
\maketitle
\section{\bf Introduction, description of the model, review of the literature and motivation}

We study a multiserver retrial queueing system having the
following structure.

\smallskip
$\bullet$ The arrival process $A(t)$ is a point process, the
increments of which form a strictly stationary and ergodic
sequence of random variables.

\smallskip
$\bullet$ There are $m$ servers, and an arriving customer occupies
one of free servers.

\smallskip
$\bullet$ If upon arrival all servers are busy, then the customer
goes to the secondary queue, orbit, and after some random time
retries more and more to occupy a server.

\smallskip
$\bullet$ A service time of each customer is exponentially
distributed random variable with parameter $\mu_1$.

\smallskip
$\bullet$ A time between retrials is exponentially distributed
with parameter $\mu_2$ for each customer in the orbit.

\smallskip
Using a martingale approach the paper provides an analysis of this
system. The paper establishes the stability condition and studies
a behavior of the limiting queue-length distributions as $\mu_2$
increases to infinity. As $\mu_2\to\infty$, the paper also proves
the convergence of appropriate queue-length distributions to those
of the associated `usual' multiserver queueing system
$A/M/m/\infty$ (without retrials), where the first parameter $A$
in the first position of the notation denotes the arrival point
process $A(t)$. In the following, by `usual' multiserver queueing
system we mean the abovementioned $A/M/m/\infty$ queueing system.
Our asymptotic results can be applied to various problems
associated with multiserver retrial queues. For example, they can
be used in performance analysis of real communication systems.

Analysis of multiserver retrial queueing systems is very hard. For
the $M/M/m$ retrial queueing systems, analytic results for the
stationary probabilities are not simple even in the case of $m=2$.
The results associated with numerical analysis have been obtained
in a large number of papers (see, e.g. Anisimov, and Artalejo
[4],
Artalejo, and Pozo
[6], Falin
[16], Neuts, and Rao
[42], Stepanov
[48], Wilkinson
[51] and others). The methods of these papers are based on
truncation of the state space for the stationary probabilities and
constructing auxiliary models helping to approximate the initial
system (see the review of Artalejo, and Falin
[5] as well as the book of Falin, and Templeton
[17] for details).

In the present paper we study a non-Markovian retrial queueing
system, the input process of which is a point process with {\it
strictly stationary and ergodic increments}. By point process with
strictly stationary and ergodic increments we mean the following.
Let $\{\xi_i\}_{i\ge 1}$ be a strictly stationary and ergodic
sequence of positive random variables, and let
$x_n=\sum_{i=1}^n\xi_i$ \ ($x_0=0$) be the corresponding sequence
of points. Then, the process $X(t)=\sum_{i=1}^\infty{\bf
I}\{x_i\le t\}$, where ${\bf I}\{\Xi\}$ denotes the indicator of
set $\Xi$, is called a point process with strictly stationary and
ergodic increments. If $\{\xi_i\}_{i\ge 1}$ is a sequence of
independent identically distributed random variables, then $X(t)$
is called a point process with independent identically distributed
increments or renewal process.

Let $\mathbb{E}\xi_1=\lambda^{-1}$. Then the assumption that
$A(t)$ is a point process with strictly stationary and ergodic
increments means that
\begin{equation}\label{1.1}
\mathbb{P}\Big\{\lim_{t\to\infty}\frac{A(t)}{t}=\lambda\Big\}=1.
\end{equation}

Along with asymptotic behavior of the limiting queue-length
distribution as parameter $\mu_2$ increases to infinity, the paper
proves continuity of the limiting queue-length distributions in
the sense of convergence of functional of the queue-length
distribution to that of the `usual' $A/M/m/\infty$ queue.

The significance of our results is motivated as follows. In the
queueing literature the multiserver $M/GI/m/0$ and $GI/M/m/0$ loss
queueing systems are the systems of special attention. Both these
queueing systems are known as a good model for telephone systems
and has been an object of investigations during many decades. The
earliest investigations of these systems were due to Palm
[43] and later due to Pollatzek
[44], Cohen
[12], Sevastyanov
[46], Tak\'acs
[49], Iglehart, and Whitt
[23]-[25] and others. Recently, the increasing attention has  been
to non-stationary multiserver loss systems (e.g. Davis, Massey,
and Whitt
[15], Massey, and Whitt
[40]) and multiserver systems with multiple customers classes
(e.g. Cohen
[13], Gail, Hantler, and Taylor
[18],
[19], Righter
[46] and others). Both these directions for multiserver queueing
system are closely related to multiserver queueing systems with
retrials and abandonments, which recently have been intensively
studied in the literature in a framework of analysis of call
centers (see Garnett, Mandelbaum, and Reiman
[21], Gans, Koole, and Mandelbaum
[20], Grier {\it et al.}
[22], Koole, and Mandelbaum
[30], Mandelbaum {\it et al.}
[36]-
[38] and others). In a framework of mentioned models of queues,
the place of the model considered in this paper is  clear. Our
assumption that input process is the process with strictly
stationary and ergodic increments is more general than those
considered earlier. Even multiserver retrial queueing models with
recurrent input are very difficult to analysis. The explicit
results for the stationary distributions of these systems are
unknown. The known results related to Markovian multiserver
retrial queueing models are not sufficient, since the real models
arising in practice not always can be good approximated by
Markovian models. Note, that Choi, Chang, and Kim
[11] applied a not standard $MAP_1, MAP_2/M/m$ retrial model to
cellular networks.

 For such general non-Markovian models as the model considered in
the paper only the stability results are an object of
investigation in the literature (see e.g. Altman, and Borovkov
[3] and references therein). The most relevant model is a model
including both retrials and abandonments, as it has been
considered in the aforementioned papers, associated with call
centers. We would like to point out that such more extended model
can also be studied by development of the method of this work. We
hope that this will be done in the future.

The analysis of the present paper is based on martingale approach.
Nowadays the martingale approach, associated with analysis of
different queueing systems and network, is familiar. Among the
well-known general textbooks on martingale theory such as  Jacod,
and Shiryayev
[26], Karatzas, and Shreve
[27], Liptser, and Shiryaev
[34],
[35], Revuz, and Yor
[45], there are special textbooks on martingale theory associated
with point processes and queues and networks such as Bremaud
[10], Baccelli, and Bremaud
[7], Whitt
[50] and others. Also there is a large number of papers studying
different queueing systems and networks with the aid of stochastic
calculus. Traditionally, the martingale methods are used to
provide weak convergence results and diffusion approximations, and
the majority of papers establish such type of results (e.g.
Abramov
[1], Kogan, and Liptser
[28], Kogan, Liptser, and Shenfild
[29], Krichagina
[31], Krichagina, Liptser, and Puhalskii
[32], Krylov, and Liptser
[33], Mandelbaum, and Pats
[39], Williams
[52] and others).

In some recent papers the martingale theory is used for analysis
of point processes and queue-length characteristics of queues and
networks also under the light traffic conditions (e.g. Abramov
[1],
[2], Kogan, and Liptser
[28], Miyazawa
[41] and references therein).

In the present paper, we study a behavior of the queue-length
process under the light traffic condition for the  multiserver
retrial queueing system, by using the known methods of the theory
of martingales. The advantage of the martingale approach is that,
it provides a deepen analysis of the system helping to study a
more wide its extension, than the traditional methods.

The paper is structured as follows. There are 11 sections. The
main results of this paper are given in Section 8. Theorem
\ref{thm8.1} establishes a property of the limiting joint
probabilities of the queue-length processes in the main queue and
orbit as parameter $\mu_2$ increases to infinity. This property is
then used in Theorem \ref{thm8.2} stating on the continuity
property of the queue-length processes, a convergence of the joint
queue-length distribution of the multiserver retrial queueing
system to that of the one-dimensional queue-length distribution of
the `usual' $A/M/m/\infty$ queueing system. Section 2 discusses
the basic equations, which are then used throughout the paper.
Section 3 deduces the Doob-Meyer semimartingale decomposition of
the basic equations. Section 4 studies normalized queue-length
processes and establishes the condition for stability (in the
sense defined precisely in Theorem \ref{thm4.1}). Section 5
derives equation for the queue-length distribution given by
Theorem \ref{thm5.1}. Section 6 is devoted to the proof of Theorem
\ref{thm5.1}. There are two corollaries of Theorem \ref{thm5.1}
given in Section 7. Sections 9 and 10 discuss algorithm for
numerical solution of the main system of equations. Specifically,
Section 10 provides numerical results under conditions of high
retrial rate, which enable us to discuss the results on
convergence obtained in Section 8. Concluding remarks are given in
Section 11.

\section{\bf Discussion of the basic equations}

All point processes considered in this paper are assumed to be
right-continuous having the left-side limits.

Consider the queue-length process of our retrial system. The
number of servers, occupied in time $t$, are denoted $Q_1(t)$, and
respectively, $Q_2(t)$ is the number of customers in orbit in time
$t$. The both queue-length processes $Q_1(t)$ and $Q_2(t)$ are
assumed to be continuous in 0, $Q_1(0)=Q_2(0)=0$, as well as
right-continuous having the left-side limits. The following two
equations describe a dynamic of the queue-length processes
$Q_1(t)$ and $Q_2(t)$. The first equation is
\begin{equation}\label{2.1}
Q_1(t)+Q_2(t)=A(t)-\int_0^t\sum_{i=1}^m{\bf I}\{Q_1(s-)\ge i\}
\mbox{d}\pi_i^{(1)}(s),
\end{equation}
where $\pi_i^{(1)}$, $i=1,2,...,m$, are independent Poisson
processes with rate $\mu_1$. The second equation is
\begin{equation}\label{2.2}
\begin{aligned}
Q_2(t)&=\int_0^t{\bf I}\{Q_1(s-)=m\}\mbox{d}A(s)\\
&-\int_0^t{\bf I}\{Q_1(s-)\neq m\}\sum_{i=1}^\infty{\bf
I}\{Q_2(s-)\ge i \}\mbox{d}\pi_i^{(2)}(s),
\end{aligned}
\end{equation}
where $\pi_i^{(2)}$, $i=1,2,...$, are independent Poisson
processes with rate $\mu_2$.

Equations \eqref{2.1} and \eqref{2.2} can be explained as follows.
The term $\sum_{i=1}^m$ ${\bf I}\{Q_1(s-)\ge i\}$ of the integrand
of \eqref{2.1} means the number of occupied servers {\it
immediately before} time $s$ in the main queue. We use the term
{\it `immediately before'} keeping in mind that the point $s$ can
be a point of possible jump. Then the right-hand side of equation
\eqref{2.1} for $Q_1(t)+Q_2(t)$ includes the number of arrivals
until time $t$ minus the number of departures given by the term
\begin{equation}\label{2.3}
\int_0^t\sum_{i=1}^m{\bf I}\{Q_1(s-)\ge i\}
\mbox{d}\pi_i^{(1)}(s).
\end{equation}
The term $\sum_{i=1}^\infty{\bf I}\{Q_2(s-)\ge i\}$ of the second
integrand of \eqref{2.2} means the number of customers in orbit
immediately before time $s$. Obviously, if there is no customer in
orbit, then the second integrand of \eqref{2.2} becomes equal to
0. Next, the term ${\bf I}\{Q_1(s-)\neq m\}$ of the second
integrand of \eqref{2.2} means that if immediately before time $s$
there is at least one free server in the main system, then one of
the customers of the orbit queue can occupy the server in time
$s$, otherwise the integrand becomes equal to 0. The first
integral of the right-hand side of \eqref{2.2} means the number of
arrivals to the orbit system during the time interval [0,$t$].
Then the right-hand side of equation \eqref{2.2} for $Q_2(t)$
includes the number of arrivals until time $t$ to the orbit minus
the number of departures from the orbit to the main queue, where
the mentioned number of arrivals to the orbit is given by the
first integral, and the mentioned number of departures from the
orbit is given by the second one.

Notice, that the equations similar to \eqref{2.1} and \eqref{2.2},
associated with time-dependent Markovian model with abandonments
and retrials, has already been considered in the literature (e.g.
Mandelbaum {\em et al.}
[36]-
[38]). However, the main emphasis of these papers was done to the
analysis of fluid limits and diffusion approximations.

\section{\bf Semimartingale decomposition of the queue-length process}

In this section we provide another representation for the
queue-length processes by using the Doob-Meyer semimartingale
decomposition (see e.g. Liptser, and Shiryayev
[34], Jacod, and Shiryayev
[26]). The compensators of the semimartingales associated with
point processes will be provided by `hat'. For example, the point
process $A(t)$ is a semimartingale, and $\widehat A(t)$ is its
compensator. The Doob-Meyer semimartingale decomposition for some
semimartingale $X$ will be written as $X=\widehat X+M_X$, where
$M_X$ is the notation for the local square-integrable martingale.
For example, the semimartingale decomposition of $A(t)$ is written
as $A(t)=\widehat A(t)+M_{A}(t)$. Along with the notation $M_X$
sometimes it is also used $M_i(t)$ or $M_{i,j}(t)$. In such cases
the sense of the local square integrable martingales $M_i(t)$ or
$M_{i,j}(t)$ is specially explained.

It is assumed in the paper that all point processes are adapted
with respect to the filtration $\mathscr{F}_t$ given on stochastic
basis $\{\Omega$, $\mathscr{F}$, ${\bf F}=(\mathscr{F}_t)_{t\ge
0}$, $\mathbb{P}\}$.

Let us start from equation \eqref{2.1}. As semimartingales, the
processes $A(t)$ and  $C(t)=\int_0^t\sum_{i=1}^m{\bf
I}\{Q_1(s-)\ge i\} \mbox{d}\pi_i^{(1)}(s)$ are represented as
\begin{equation}\label{3.1}
A(t)=\widehat A(t)+M_A(t),
\end{equation}
\begin{equation}\label{3.2}
C(t)=\widehat C(t)+M_C(t),
\end{equation}
The compensator $\widehat C(t)$ has the representation
\begin{equation}\label{3.3}
\widehat C(t)=\mu_1\int_0^tQ_1(s) \mbox{d}s
\end{equation}
(for details see Dellacherie
[14], Liptser, and Shiryayev
[34],
[35], Theorem 1.6.1).

By virtue of \eqref{3.1}, \eqref{3.2}, and \eqref{3.3}, for
equation \eqref{2.1}, we have
\begin{equation}\label{3.4}
Q_1(t)+Q_2(t)=\widehat A(t)+M_A(t)-\mu_1\int_0^tQ_1(s)
\mbox{d}s-M_C(t).
\end{equation}
Denoting $M_1(t)=M_A(t)-M_C(t)$ we obtain
\begin{equation}\label{3.5}
Q_1(t)+Q_2(t)=\widehat A(t)-\mu_1\int_0^tQ_1(s) \mbox{d}s+M_1(t).
\end{equation}

Let us now consider equation \eqref{2.2}. For the associated
arrival process $D_1(t)$ = $\int_0^t{\bf
I}\{Q_1(s-)=m\}\mbox{d}A(s)$ we have the following:
\begin{equation}\label{3.6}
D_1(t)=\widehat D_1(t)+M_{D_1}(t),
\end{equation}
where
\begin{equation}\label{3.7}
\widehat D_1(t)=\int_0^t{\bf I}\{Q_1(s-)=m\}\mbox{d}\widehat A(s)
\end{equation}
and
\begin{equation}\label{3.8}
M_{D_1}(t)=\int_0^t{\bf I}\{Q_1(s-)=m\}\mbox{d}M_A(s).
\end{equation}

For the associated departure process $D_2(t)$ = $\int_0^t{\bf
I}\{Q_1(s-)\neq m\}$ $\sum_{i=1}^\infty$ ${\bf I}\{Q_2(s-)$ $\ge i
\}\mbox{d}\pi_i^{(2)}(s)$ we have the following:
\begin{equation}\label{3.9}
D_2(t)=\widehat D_2(t)+M_{D_2}(t),
\end{equation}
where
\begin{equation}\label{3.10}
\widehat D_2(t)=\mu_2\int_0^t{\bf I}\{Q_1(s)\neq
m\}Q_2(s)\mbox{d}s
\end{equation}
(see Dellacherie
[14], Liptser, and Shiryayev
[34],
[35], Theorem 1.6.1).

Then \eqref{2.2} can be rewritten as follows:
\begin{equation}\label{3.11}
\begin{aligned}
Q_2(t)&=\int_0^t{\bf I}\{Q_1(s-)=m\}\mbox{d}\widehat A(s)\\
&- \mu_2\int_0^t{\bf I}\{Q_1(s)\neq m\}Q_2(s)\mbox{d}s +M_2(t),
\end{aligned}
\end{equation}
where
\begin{equation}\label{3.12}
M_2(t)=M_{D_1}(t)-M_{D_2}(t).
\end{equation}

\section{\bf Normalized queue-length processes and condition for the stability}

In this section we study the normalized queue-length processes
\begin{equation}\label{4.1}
q_k(t)=\frac{Q_k(t)}{t}, \ k=1,2; \ t>0,
\end{equation}
and its asymptotic properties as $t\to\infty$. In the following
the small letters will stand for normalized processes. The
notation for normalized processes corresponds to the notation of
original processes given by capital letters. For example,
$\widehat a(t)=t^{-1}$ $\widehat A(t)$; $m_{D_2}(t)=t^{-1}$ $
M_{D_2}(t)$ and so on.

Following this comment, the equations associated with the
queue-length processes \eqref{3.5} and \eqref{3.11} can be written
\begin{equation}\label{4.2}
q_1(t)+q_2(t)=\widehat a(t)-\frac{\mu_1}{t}\int_0^tsq_1(s)
\mbox{d}s+m_1(t),
\end{equation}
and
\begin{equation}\label{4.3}
\begin{aligned}
q_2(t)&=\frac{1}{t}\int_0^t{\bf I}\{Q_1(s-)=m\}\mbox{d}[s\widehat
a(s)]\\
&- \frac{\mu_2}{t}\int_0^t{\bf I}\{Q_1(s)\neq m\}sq_2(s)\mbox{d}s
+m_2(t),
\end{aligned}
\end{equation}

\smallskip
Let us now study these two equations \eqref{4.2} and \eqref{4.3}
as $t\to\infty$. More accurately, let us find $\mathbb{P}\lim$ of
$q_1(t)$ and $q_2(t)$ as $t\to\infty$, where $\mathbb{P}\lim$
denotes the limit in probability.

Show that
\begin{equation}\label{4.4}
\mathbb{P}\lim_{t\to\infty}m_1(t)=0.
\end{equation}
Indeed, because of $m_1(t)=m_A(t)-m_C(t)$, we have
\begin{equation}\label{4.5}
\mathbb{P}\lim_{t\to\infty}m_1(t)\le \mathbb{
P}\lim_{t\to\infty}\big|m_A(t)\big|+\mathbb{
P}\lim_{t\to\infty}\big|m_C(t)\big|.
\end{equation}
Applying the Lenglart-Rebolledo inequality we obtain:
\begin{equation}\label{4.6}
\begin{aligned}
 \mathbb{P}\{|m_A(t)|>\delta\}&\leq\mathbb{P}\Big\{\sup_{0<s\leq
t}\Big|\frac{sm_A(s)}{t}\Big|>\delta \Big\}\\
&=\mathbb{P}\Big\{\sup_{0<s\le t}\big|M_A(s)\big|>\delta t
\Big\}\\
&\leq\frac{\epsilon}{\delta^2}+\mathbb{P}\{A(t)>\epsilon t^2\}\\
&
=\frac{\epsilon}{\delta^2}+\mathbb{P}\Big\{\frac{A(t)}{t}>\epsilon
t\Big\}
\end{aligned}
\end{equation}
The both terms of the right-hand side vanish, as $\epsilon$ is
taken sufficiently small and $t$ increases to infinity such that
$\epsilon t\to\infty$. That is $\mathbb{P}\lim_{t\to\infty}$ $
\big|m_A(t)\big|=0$.

Taking into account that
\begin{equation}\label{4.7}
\int_0^t\sum_{i=1}^m{\bf I}\{Q_1(s-)\ge i\}
\mbox{d}\pi_i^{(1)}(s)\le \sum_{i=1}^m\pi_i^{(1)}(t),
\end{equation}
by virtue of the Lenglart-Rebolledo inequality we have:
\begin{equation}\label{4.8}
\begin{aligned}
\mathbb{P}\{|m_C(t)|>\delta\}&\leq\mathbb{P}\Big\{\sup_{0<s\le
t}\Big|\frac{sm_C(s)}{t}\Big|>\delta \Big\}\\
& =\mathbb{P}\Big\{\sup_{0<s\le t}\big|M_C(s)\big|>\delta t
\Big\}\\& \leq\frac{\epsilon}{\delta^2}+\mathbb{
P}\Big\{\sum_{i=1}^m\pi_{i}^{(1)}(t)>\epsilon t^2\Big\}\\
&= \frac{\epsilon}{\delta^2}+\mathbb{
P}\Big\{\frac{1}{t}\sum_{i=1}^m\pi_{i}^{(1)}(t)>\epsilon t\Big\}
.\end{aligned}
\end{equation}
 As earlier (see reference \eqref{4.6}), now we obtain
$\mathbb{P}\lim_{t\to\infty}\big|m_C(t)\big|=0$. Thus, it is shown
that $\mathbb{P}\lim_{t\to\infty}m_1(t)=0$.

Analogously to the above, $m_2(t)=m_{D_1}(t)-m_{D_2}(t)$.
Therefore, similarly to \eqref{3.5}
\begin{equation}\label{4.9}
\mathbb{P}\lim_{t\to\infty}m_2(t)\le \mathbb{
P}\lim_{t\to\infty}\big|m_{D_1}(t)\big|+\mathbb{
P}\lim_{t\to\infty}\big|m_{D_2}(t)\big|.
\end{equation}
Notice (see \eqref{3.8}), that $|m_{D_1}(t)|\leq|m_A(t)|$ for all
$t>0$. Therefore,
\begin{equation}\label{4.10}
\mathbb{P}\lim_{t\to\infty}\big|m_{D_1}(t)\big|\le\mathbb{
P}\lim_{t\to\infty}\big|m_A(t)\big|=0.
 \end{equation}
 However,
both $\mathbb{P}\lim_{t\to\infty}|m_{D_2}(t)|=0$ and
$\mathbb{P}\lim_{t\to\infty}|m_2(t)|=0$ can only be true under the
additional condition $\mathbb{P}\lim_{t\to\infty}a(t)<\mu_1m$.
Recall that according to \eqref{1.1} $\mathbb{
P}\lim_{t\to\infty}a(t)=\lambda$. Therefore, the above condition
is $\lambda<\mu_1m$.

In order to prove $\mathbb{P}\lim_{t\to\infty}|m_2(t)|=0$ and
$\mathbb{P}\lim_{t\to\infty}|m_{D_2}(t)|=0$ under the
abovementioned additional condition $\lambda<\mu_1m$ let us now
study equation \eqref{4.2}.

Notice first that from the fact that $\mathbb{
P}\lim_{t\to\infty}|m_A(t)|=0$ we have
\begin{equation}\label{4.11}
\mathbb{P}\lim_{t\to\infty}\widehat a(t)=\mathbb{
P}\lim_{t\to\infty}a(t)=\lambda,
\end{equation}
since
\begin{equation}\label{4.12}
\mathbb{P}\lim_{t\to\infty}\widehat a(t)=\mathbb{
P}\lim_{t\to\infty}a(t)-\mathbb{P}\lim_{t\to\infty}m_A(t).
\end{equation}
Next, if $\lim_{t\to\infty}t^{-1}\int_0^t\mathbb{
P}\{Q_1(s)=m\}\mbox{d}s=1$, then by the Lebesgue dominated
convergence
\begin{equation}\label{4.13}
\mathbb{E}\lim_{t\to\infty}\frac{1}{t}\int_0^t{\bf
I}\{Q_1(s)=m\}\mbox{d}[s\widehat a(s)]=\lambda,
\end{equation}
and
\begin{equation}\label{4.14}
\lim_{t\to\infty}\mathbb{E}q_2(t)=\lim_{t\to\infty}\mathbb{E}\widehat
a(t)=\lambda.
\end{equation}
\eqref{4.14} means that if
$\lim_{t\to\infty}t^{-1}\int_0^t\mathbb{
P}\{Q_1(s)=m\}\mbox{d}s=1$, then $\mathbb{E} Q_2(t)$ increases to
infinity, as $t\to\infty$.

Let us now assume, that $\lim_{t\to\infty}t^{-1}\int_0^t\mathbb{
P}\{Q_1(s)=m\}\mbox{d}s=p<1.$ Then, it is not difficult to prove
that only $\mathbb{ P}\lim_{t\to\infty}q_2(t)=0$ must satisfy
\eqref{4.2}.

Indeed, assume that $\mathbb{P}\lim_{t\to\infty}q_2(t)>0$. Then,
taking into account that $\lim_{t\to\infty}$ $
t^{-1}\int_0^t\mathbb{ P}\{Q_1(s)\neq m\}\mbox{d}s=1-p>0$, for
large $t$ we obtain
\begin{equation}\label{4.15}
\int_0^t{\bf
I}\{Q_1(s)\neq m\}sq_2(s)\mbox{d}s=O(t^2).
\end{equation}

This means that
\begin{equation}\label{4.16}
\mathbb{P}\lim_{t\to\infty}\frac{\mu_2}{t}\int_0^t{\bf
I}\{Q_1(s)\neq m\}sq_2(s)\mbox{d}s=\infty,
\end{equation}
 and
therefore $\lim_{t\to\infty}\mathbb{E}q_2(t)=-\infty$. This
contradicts to the fact that $\mathbb{E}q_2(t)\ge 0$ for all $t\ge
0$, and therefore, only $\mathbb{P}\lim_{t\to\infty}q_2(s)=0$ is a
possible limit in probability, satisfying \eqref{4.3} for some
unique value $
p=p^*=\lim_{t\to\infty}t^{-1}\int_0^t\mathbb{P}\{Q_1(s)=
m\}\mbox{d}s$.

Thus, we proved that if $\lim_{t\to\infty}t^{-1}\int_0^t\mathbb{
P}\{Q_1(s)=m\}\mbox{d}s=p<1$, then $\mathbb{P}\lim_{t\to\infty}$
$q_2(s)=0$. In fact, taking into account that $Q_2(t)$ is a
c\'adl\'ag process having with probability 1 a finite number of
jumps in any finite interval, from the above contradiction we have
\begin{equation}\label{4.17}
\lim_{t\to\infty}\frac{1}{t}\int_0^t\mathbb{E}Q_2(s)\mbox{d}s<\infty,
\end{equation}
and therefore,
\begin{equation}\label{4.18}
\mathbb{P}\lim_{t\to\infty}\frac{1}{t}\int_0^t{\bf
I}\{Q_2(s)<\infty\}\mbox{d}s=\lim_{t\to\infty}\frac{1}{t}\int_0^t\mathbb{
P}\{Q_2(s)<\infty\}\mbox{d}s=1.
\end{equation}

Now, using the Lenglart-Rebolledo inequality for large $t$ we
obtain

\begin{equation}\label{4.19}
\begin{aligned}
 \mathbb{P}\{|m_{D_2}(t)|>\delta\}&\leq\mathbb{P}\Big\{\sup_{o<s\le
t}\Big|\frac{sm_{D_2}(s)}{t}\Big|>\delta\Big\}\\
& =\mathbb{P}\{\sup_{o<s\le t}\big|M_{D_2}(s)\big|>\delta t\}\\
&\leq \frac{\epsilon}{\delta^2}+\mathbb{P}\Big\{\int_0^t{\bf
I}\{Q_1(s)\neq m\} sq(s)\mbox{d}s>\epsilon t^2\Big\}\\
&\le \frac{\epsilon}{\delta^2}+\mathbb{P}\Big\{\int_0^t
sq(s)\mbox{d}s>\epsilon
t^2\Big\}\\&=\frac{\epsilon}{\delta^2}+\mathbb{P}\Big\{\frac{1}{t}\int_0^t
sq(s)\mbox{d}s>\epsilon t\Big\}.
\end{aligned}
\end{equation}
 Therefore, under
the assumption that $\lim_{t\to\infty}
t^{-1}\int_0^t\mathbb{P}\{Q_1(s)=m\}\mbox{d}s$ $=p<1$, we obtain
\begin{equation}\label{4.20}
\mathbb{P}\lim_{t\to\infty}\big|m_{D_2}(t)\big|=0.
\end{equation}
Hence, together with \eqref{4.10}, under the assumption that
$\lim_{t\to\infty}$ $t^{-1}$ $\int_0^t$
$\mathbb{P}\{Q_1(s)=m\}\mbox{d}s=p<1$ we have
\begin{equation}\label{4.21}
\mathbb{P}\lim_{t\to\infty}\big|m_{2}(t)\big|=0.
\end{equation}

Let us return to relation \eqref{4.18} under the condition
$\lambda<\mu_1m$. We have
\begin{equation}\label{4.22}
\lim_{t\to\infty}\frac{1}{t}\int_0^t\sum_{l=0}^\infty\mathbb{
P}\{Q_2(s)=l\}\mbox{d}s=1,
\end{equation}
and, since $Q_1(s)$ can take values $0,1,\ldots,m$ only, we have
\begin{equation}\label{4.23}
\lim_{t\to\infty}\frac{1}{t}\int_0^t\sum_{l=0}^m\mathbb{
P}\{Q_1(s)=l\}\mbox{d}s=1.
\end{equation}
Therefore, under the condition $\lambda<\mu_1m$,
\begin{equation}\label{4.24}
\begin{aligned}
&\sum_{l=0}^m\sum_{k=0}^\infty\lim_{t\to\infty}\frac{1}{t}\int_0^t\mathbb{
P}\{Q_1(s)=l, Q_2(s)=k\}\mbox{d}s\\
&
=\lim_{t\to\infty}\frac{1}{t}\int_0^t\sum_{l=0}^m\sum_{k=0}^\infty\mathbb{
P}\{Q_1(s)=l, Q_2(s)=k\}\mbox{d}s\\
&=1.
\end{aligned}
\end{equation}
Thus, we have the following theorem.

\begin{thm}\label{thm4.1}
 Under the condition $\lambda<\mu_1m$
there exist
$$
\lim_{t\to\infty}\frac{1}{t}\int_0^t\mathbb{P}\{Q_1(s)=l,
Q_2(s)=k\}\mbox{d}s,
$$
$$
l=0,1,\ldots,m; \ k=0,1,\ldots,
$$
satisfying \eqref{4.24}.
\end{thm}

\section{\bf Analysis of the limiting queue-length distributions}

In the rest of the paper it is assumed that condition
$\lambda<\mu_1m$ is fulfilled, and therefore the system is stable
in the sense of Theorem \ref{thm4.1}. Notice, that the statement
of Theorem \ref{thm4.1} does not mean existence of the limiting
stationary probabilities as $t\to\infty$. For example, if
increments of the point process $A(t)$ are lattice, then the
stationary probabilities do not exist.

In this case one can speak only about appropriate fractions of
time in two-dimensional states ($i,j$) associated with the
queue-length processes $Q_1(t)$ and $Q_2(t)$. In general, we speak
about
\begin{equation}\label{5.1}
\lim_{t\to\infty}\frac{1}{t}\int_0^t\mathbb{P}\{Q_1(s)=i,
Q_2(s)=j\}\mbox{d}s
\end{equation}
rather than $\lim_{t\to\infty}\mathbb{P}\{Q_1(t)=i, Q_2(t)=j\}$.
However, if the increments of the point process $A(t)$ are
independent, identically distributed and non-lattice, then there
exist the limiting stationary probabilities
$\lim_{t\to\infty}\mathbb{P}\{Q_1(t)=i, ~Q_2(t)=j\}$ coinciding
with \eqref{5.1}. Indeed, in this case the process
$Q(t)=Q_1(t)+Q_2(t)$ has a structure of regeneration process, and
then the proof of this fact follows by a slight extension of
arguments given in the proofs of Theorem 5 on p. 173 and Theorem
22 on p. 157 of Borovkov
[9].

\smallskip
Let us introduce the processes
\begin{equation}\label{5.2}
I_{i,j}(t)={\bf I}\{Q_1(t)=i~\cap~Q_2(t)=j\}, ~i=0,1,...,m; \
~j=0,1,...~,
\end{equation}
assuming that $I_{-1,j}(t)\equiv 0$ and $I_{i,-1}(t)\equiv 0$.

The jump of a point process is denoted by adding $\Delta$. For
example, $\Delta A(t)$ is a jump of $A(t)$, $\Delta
\pi_k^{(1)}(t)$ is a jump of the $k$th Poisson process with rate
$\mu_1$ etc.

Let us denote:
\begin{equation}\label{5.3}
 \Pi_i^{(1)}(t)=\sum_{k=1}^i\pi_k^{(1)}(t),
\end{equation}
and
\begin{equation}\label{5.4}
\Pi_j^{(2)}(t)=\sum_{k=1}^j\pi_k^{(2)}(t).
\end{equation}
 Taking
into account that the jumps of all the processes $A(t)$,
$\pi_k^{(1)}(t)$ and $\pi_l^{(2)}(t)$ are disjoint ($k=1,2,...,m$;
$l=1,2,...$), we have the following equations:
\begin{equation}\label{5.5}
\begin{aligned}
& {\bf I}\{Q_1(t-)+\Delta Q_1(t)=i~\cap~Q_2(t)=Q_2(t-)=j\}\\
& =I_{i-1,j}(t-)\Delta A(t)
+I_{i+1,j}(t-)\Delta\Pi_{i+1}^{(1)}(t)\\
& +I_{i,j}(t-)[1-\Delta A(t)][1-\Delta \Pi_i^{(1)}(t)][1-\Delta
\Pi_j^{(2)}(t)],\\
& i=0,1,...,m-1;
\end{aligned}
\end{equation}

\medskip
\begin{equation}\label{5.6}
\begin{aligned}
& {\bf I}\{Q_1(t-)+\Delta Q_1(t)=i~\cap~Q_2(t-)\ne Q_2(t)=j\}\\
& =I_{i-1,j+1}(t-)\Delta \Pi_{j+1}^{(2)}(t),\\&i=0,1,...,m-1;
\end{aligned}
\end{equation}
\medskip

\begin{equation}\label{5.7}
\begin{aligned}
&{\bf I}\{Q_1(t-)+\Delta Q_1(t)=m~\cap~Q_2(t)=Q_2(t-)=j\}\\
& =I_{m-1,j}(t-)\Delta A(t)\\& +I_{m,j}(t-)[1-\Delta
A(t)][1-\Delta \Pi_m^{(1)}(t)];
\end{aligned}
\end{equation}
\medskip

\begin{equation}\label{5.8}
\begin{aligned}
&{\bf I}\{Q_1(t-)+\Delta Q_1(t)=m~\cap~Q_2(t-)\ne Q_2(t)=j\}\\&
=I_{m,j-1}(t-)\Delta A(t).
\end{aligned}
\end{equation}
\medskip

Then,
\begin{equation}\label{5.9}
\begin{aligned}
&\Delta I_{i,j}(t)={\bf I}\{Q_1(t-)+\Delta
Q_1(t)=i~\cap~Q_2(t)=Q_2(t-)=j\}\\ &+ {\bf I}\{Q_1(t-)+\Delta
Q_1(t)=i~\cap~Q_2(t-)\ne Q_2(t)=j\}\\ & -I_{i,j}(t-),\\
&i=0,1,...,m; \ j\ge 0.
\end{aligned}
\end{equation}
\medskip
Since
\begin{equation}\label{5.10}
\sum_{s\le t}\Delta I_{i,j}(s)=I_{i,j}(t)-I_{i,j}(0),
\end{equation}
we have the following.

For $i=0,1,...,m-1$,
\begin{equation}\label{5.11}
\begin{aligned}
I_{i,j}(t)&=I_{i,j}(0)+\int_0^t\big[I_{i-1,j}(s-)-I_{i,j}(s-)\big]\mbox{d}A(s)\\
&-
\int_0^tI_{i,j}(s-)\mbox{d}\Pi_i^{(1)}(s)-\int_0^tI_{i,j}(s-)\mbox{d}\Pi_j^{(2)}(s)\\
&+ \int_0^tI_{i+1,j}(s-)\mbox{d}\Pi_{i+1}^{(1)}(s)\\
&+\int_0^tI_{i-1,j+1}(s-)\mbox{d}\Pi_{j+1}^{(2)}(s).
\end{aligned}
\end{equation}
In turn, for $i=m$ we have
\begin{equation}\label{5.12}
\begin{aligned}
I_{m,j}(t)&=I_{m,j}(0)+\int_0^t\big[I_{m-1,j}(s-)-I_{m,j}(s-)\big]\mbox{d}A(s)\\
&-
\int_0^tI_{m,j}(s-)\mbox{d}\Pi_m^{(1)}(s)+\int_0^tI_{m,j-1}(s-)\mbox{d}A(s)\\
&+ \int_0^tI_{m-1,j+1}(s-)\mbox{d}\Pi_{j+1}^{(2)}(s).
\end{aligned}
\end{equation}

Using the Doob-Meyer semimartingale decomposition, from
\eqref{5.11} and \eqref{5.12} we obtain the following equations.
For $i=0,1,...,m-1$
\begin{equation}\label{5.13}
\begin{aligned}
I_{i,j}(t)&=I_{i,j}(0)+\int_0^t\big[I_{i-1,j}(s-)-I_{i,j}(s-)\big]\mbox{d}\widehat
A(s)\\&- \mu_1i\int_0^tI_{i,j}(s)\mbox{d}s
-\mu_2j\int_0^tI_{i,j}(s)\mbox{d}s\\&+\mu_1(i+1)\int_0^tI_{i+1,j}(s)\mbox{d}s\\&+
\mu_2(j+1)\int_0^tI_{i-1,j+1}(s)\mbox{d}s+M_{i,j}(t),
\end{aligned}
\end{equation}
where the local square integrable martingale $M_{i,j}(t)$ has the
representation
\begin{equation}\label{5.14}
\begin{aligned}
M_{i,j}(t)&=\int_0^t\big[I_{i-1,j}(s-)-I_{i,j}(s-)\big]\mbox{d}[A(s)-\widehat
A(s)]\\&- \int_0^tI_{i,j}(s-)\mbox{d}[\Pi_i^{(1)}(s)-\mu_1is]\\&
-\int_0^tI_{i,j}(s-)\mbox{d}[\Pi_j^{(1)}(s)-\mu_2js]\\&
+\int_0^tI_{i+1,j}(s-)\mbox{d}[\Pi_{i+1}^{(1)}(s)-(i+1)\mu_1s]\\&+
\int_0^tI_{i-1,j+1}(s-)\mbox{d}[\Pi_{j+1}^{(2)}(s)-(j+1)\mu_2s].
\end{aligned}
\end{equation}
In turn, for $i=m$ we have
\begin{equation}\label{5.15}
\begin{aligned}
I_{m,j}(t)&=I_{m,j}(0)+\int_0^t\big[I_{m-1,j}(s-)-I_{m,j}(s-)\big]\mbox{d}\widehat
A(s)\\ &-
\mu_1m\int_0^tI_{m,j}(s)\mbox{d}s+\int_0^tI_{m,j-1}(s-)\mbox{d}\widehat
A(s)\\&+ \mu_2(j+1)\int_0^tI_{m-1,j+1}(s)\mbox{d}s+M_{m,j}(t),
\end{aligned}
\end{equation}
where the local square integrable martingale $M_{m,j}(t)$ has the
representation
\begin{equation}\label{5.16}
\begin{aligned}
M_{m,j}(t)&=\int_0^t\big[I_{m-1,j}(s-)-I_{m,j}(s-)\big]\mbox{d}[A(s)-\widehat
A(s)]\\&- \int_0^tI_{m,j}(s-)\mbox{d}[\Pi_m^{(1)}(s)-\mu_1ms]\\
&+\int_0^tI_{m,j-1}(s-)\mbox{d}[A(s)-\widehat A(s)]\\
&+
\int_0^tI_{m-1,j+1}(s-)\mbox{d}[\Pi_{j+1}^{(2)}(s)-\mu_2(j+1)s].
\end{aligned}
\end{equation}
Now, we are ready to formulate and prove the theorem.

\begin{thm}\label{thm5.1}
For the queue-length processes $Q_1(t)$ and
$Q_2(t)$ we have the following equations.\\

(i) In the case $i=0$
\begin{equation}\label{5.17}
\begin{aligned}
&\lim_{t\to\infty}\frac{1}{t}\mathbb{E}\int_0^t{\bf I}\{Q_1(s-)=0,
Q_2(s-)=j\}\mbox{d}A(s)\\
& =\mu_1 \lim_{t\to\infty}\frac{1}{t}\int_0^t\mathbb{P}\{Q_1(s)=1,
Q_2(s)=j\}\mbox{d}s\\&
-\mu_2j\lim_{t\to\infty}\frac{1}{t}\int_0^t\mathbb{P}\{Q_1(s)=0,
Q_2(s)=j\}\mbox{d}s.
\end{aligned}
\end{equation}

(ii) In the case $i=1,2,...,m-1$ ($m\ge 2$)
\begin{equation}\label{5.18}
\begin{aligned}
&\lim_{t\to\infty}\frac{1}{t}\mathbb{E}\int_0^t[{\bf
I}\{Q_1(s-)=i, Q_2(s-)=j\}\\& -{\bf I}\{Q_1(s-)=i-1,
Q_2(s-)=j\}]\mbox{d}A(s)\\&=
\mu_1(i+1)\lim_{t\to\infty}\frac{1}{t}\int_0^t\mathbb{P}\{Q_1(s)=i+1,
Q_2(s)=j\}\mbox{d}s\\&
-\mu_1i\lim_{t\to\infty}\frac{1}{t}\int_0^t\mathbb{P}\{Q_1(s)=i,
Q_2(s)=j\}\mbox{d}s\\&
-\mu_2j\lim_{t\to\infty}\frac{1}{t}\int_0^t\mathbb{P}\{Q_1(s)=i,
Q_2(s)=j\}\mbox{d}s\\&
+\mu_2(j+1)\lim_{t\to\infty}\frac{1}{t}\int_0^t\mathbb{P}\{Q_1(s)=i-1,
Q_2(s)=j+1\}\mbox{d}s.
\end{aligned}
\end{equation}

(iii) In the case $i=m$
\begin{equation}\label{5.19}
\begin{aligned}
&\lim_{t\to\infty}\frac{1}{t}\mathbb{E}\int_0^t[{\bf
I}\{Q_1(s-)=m,
Q_2(s-)=j\}\\& -{\bf I}\{Q_1(s-)=m-1, Q_2(s-)=j\}\\
&-{\bf I}\{Q_1(s-)=m, Q_2(s-)=j-1\} ]\mbox{d}A(s)\\
&=
\mu_2(j+1)\lim_{t\to\infty}\frac{1}{t}\int_0^t\mathbb{P}\{Q_1(s)=m-1,
Q_2(s)=j+1\}\mbox{d}s\\&
-\mu_1m\lim_{t\to\infty}\frac{1}{t}\int_0^t\mathbb{P}\{Q_1(s)=m,
Q_2(s)=j\}.
\end{aligned}
\end{equation}
Here in \eqref{5.17}-\eqref{5.19} it is put ${\bf I}\{Q_1(t)=i,
Q_2(t)=j\}=0$ if at least one of the values $i$, $j$ is equal to
-1.
\end{thm}

\section{\bf Proof of Theorem \ref{thm5.1}}

Let us first study the case $i=0$. From \eqref{5.13} we have
\begin{equation}\label{6.1}
\begin{aligned}
\mathbb{P}
\lim_{t\to\infty}\frac{1}{t}\Big(I_{0,j}(t)-I_{0,j}(0)\Big)=&-\mathbb{
P}\lim_{t\to\infty}\frac{1}{t}\int_0^tI_{0,j}(s-)\mbox{d}\widehat
A(s)\\
& +\mathbb{P}
\lim_{t\to\infty}\frac{1}{t}\mu_1\int_0^tI_{1,j}(s)\mbox{d}s\\
&-\mathbb{P}
\lim_{t\to\infty}\frac{1}{t}\mu_2j\int_0^tI_{0,j}(s)\mbox{d}s\\
& +\mathbb{P}\lim_{t\to\infty}\frac{M_{0,j}(t)}{t}.
\end{aligned}
\end{equation}
The left-hand side of \eqref{6.1} is equal to zero. Therefore,
rewriting \eqref{6.1} in the form $0=-K_1+K_2-K_3+K_4$, let us
compute the terms of the right-hand side. Using the Lebesgue
theorem on dominated convergence we have
\begin{equation}\label{6.2}
\begin{aligned}
K_1&=\mathbb{P}
\lim_{t\to\infty}\frac{1}{t}\int_0^tI_{0,j}(s-)\mbox{d}\widehat
A(s)\\& = \lim_{t\to\infty}\frac{1}{t}\mathbb{
E}\int_0^tI_{0,j}(s-)\mbox{d}\widehat A(s).
\end{aligned}
\end{equation}
Taking into account that $A(t)/t$ and $\widehat A(t)/t$ have the
same limit in probability (see  \eqref{4.11}), relation
\eqref{6.2} can be finally rewritten as follows:
\begin{equation}\label{6.3}
\begin{aligned}
K_1&=\lim_{t\to\infty}\frac{1}{t}\mathbb{E}\int_0^t{\bf
I}\{Q_1(s-)=0, Q_2(s-)=j\}\mbox{d}\widehat A(s)\\ &=
\lim_{t\to\infty}\frac{1}{t}\mathbb{E}\int_0^t{\bf I}\{Q_1(s-)=0,
Q_2(s-)=j\}\mbox{d}[A(s)-M_A(s)]\\&=
\lim_{t\to\infty}\frac{1}{t}\mathbb{E}\int_0^t{\bf I}\{Q_1(s-)=0,
Q_2(s-)=j\}\mbox{d}A(s).
\end{aligned}
\end{equation}

Next,
\begin{equation}\label{6.4}
\begin{aligned}
K_2&=\mathbb{P}
\lim_{t\to\infty}\frac{1}{t}\mu_1\int_0^tI_{1,j}(s)\mbox{d}s\\&=
\mu_1\lim_{t\to\infty}\frac{1}{t}\int_0^t\mathbb{P}\{Q_1(s)=1,
Q_2(s)=j\}\mbox{d}s.
\end{aligned}
\end{equation}
Similarly,
\begin{equation}\label{6.5}
\begin{aligned}
K_3&=\mathbb{P}-
\lim_{t\to\infty}\frac{1}{t}\mu_2j\int_0^tI_{0,j}(s)\mbox{d}s\\
&= \mu_2j\lim_{t\to\infty}\frac{1}{t}\int_0^t\mathbb{P}\{Q_1(s)=0,
Q_2(s)=j\}\mbox{d}s.
\end{aligned}
\end{equation}
Notice, that if $j=0$ then $K_2=0$. Next,
\begin{equation}\label{6.6}
\begin{aligned}
K_4&= \mathbb{P}\lim_{t\to\infty}\Big|\frac{M_{i,j}(t)}{t}\Big|\\
& \leq \mathbb{
P}\lim_{t\to\infty}\Big(|m_A(t)|+\Big|\frac{\pi_1^{(1)}(t)-\mu_1t}{t}\Big|\\&+
\Big|\frac{\Pi_j^{(2)}(t)-\mu_2jt}{t}\Big|
+\Big|\frac{\Pi_{j+1}^{(2)}(t)-\mu_2(j+1)t}{t}\Big|\Big)\\&=0.
\end{aligned}
\end{equation}

Thus, for $i=0$ we have the following
\begin{equation}\label{6.7}
\begin{aligned}
&\lim_{t\to\infty}\frac{1}{t}\mathbb{E}\int_0^t{\bf I}\{Q_1(s-)=0,
Q_2(s-)=j\}\mbox{d}A(s)\\& =\mu_1
\lim_{t\to\infty}\frac{1}{t}\int_0^t\mathbb{P}\{Q_1(s)=1,
Q_2(s)=j\}\mbox{d}s\\&
-\mu_2j\lim_{t\to\infty}\frac{1}{t}\int_0^t\mathbb{P}\{Q_1(s)=0,
Q_2(s)=j\}\mbox{d}s.
\end{aligned}
\end{equation}
\eqref{5.17} follows.

\smallskip
Let us consider now the case $1\le i\le m-1$ ~($m\ge 2$). We have
the following equation:
\begin{equation}\label{6.8}
\begin{aligned}
\mathbb{P}
\lim_{t\to\infty}\frac{1}{t}\Big(I_{i,j}(t)-I_{i,j}(0)\Big)
&=\mathbb{P}
\lim_{t\to\infty}\frac{1}{t}\int_0^tI_{i-1,j}(s-)\mbox{d}\widehat
A(s)\\&-\mathbb{P}
\lim_{t\to\infty}\frac{1}{t}\int_0^tI_{i,j}(s-)\mbox{d}\widehat
A(s)\\& +\mathbb{P}
\lim_{t\to\infty}\frac{1}{t}\mu_1(i+1)\int_0^tI_{i+1,j}(s)\mbox{d}s\\
&-\mathbb{P}
\lim_{t\to\infty}\frac{1}{t}\mu_2j\int_0^tI_{i,j}(s)\mbox{d}s\\
&-\mathbb{
P}\lim_{t\to\infty}\frac{1}{t}\mu_1i\int_0^tI_{i,j}(s)\mbox{d}s\\
&+\mathbb{
P}\lim_{t\to\infty}\frac{1}{t}\mu_2(j+1)\int_0^tI_{i-1,j+1}(s)\mbox{d}s\\
&+\mathbb{P}\lim_{t\to\infty}\frac{M_{i,j}(t)}{t},
\end{aligned}
\end{equation}
The term of the left-hand side of \eqref{6.8} is equal to zero.
Therefore rewriting \eqref{6.8} as $0=K_1-K_2+K_3-K_4-K_5+K_6+K_7$
we have the following. Similarly to \eqref{6.3}
\begin{equation}\label{6.9}
K_1=\lim_{t\to\infty}\frac{1}{t}\mathbb{E}\int_0^t{\bf
I}\{Q_1(s-)=i-1, Q_2(s-)=j\}\mbox{d}A(s),
\end{equation}
and
\begin{equation}\label{6.10}
K_2=\lim_{t\to\infty}\frac{1}{t}\mathbb{E}\int_0^t{\bf
I}\{Q_1(s-)=i, Q_2(s-)=j\}\mbox{d}A(s).
\end{equation}
Similarly to \eqref{6.4} and \eqref{6.5}
\begin{equation}\label{6.11}
K_3=\mu_1(i+1)\lim_{t\to\infty}\frac{1}{t}\int_0^t\mathbb{
P}\{Q_1(s)=i+1, Q_2(s)=j\}\mbox{d}s,
\end{equation}
\begin{equation}\label{6.12}
K_4=\mu_2j\lim_{t\to\infty}\frac{1}{t}\int_0^t\mathbb{P}\{Q_1(s)=i,
Q_2(s)=j\}\mbox{d}s,
\end{equation}
\begin{equation}\label{6.13}
K_5=\mu_1i\lim_{t\to\infty}\frac{1}{t}\int_0^t\mathbb{P}\{Q_1(s)=i,
Q_2(s)=j\}\mbox{d}s,
\end{equation}
\begin{equation}\label{6.14}
K_6=\mu_2(j+1)\lim_{t\to\infty}\frac{1}{t}\int_0^t\mathbb{
P}\{Q_1(s)=i-1, Q_2(s)=j+1\}\mbox{d}s,
\end{equation}
and similarly to \eqref{6.6}
\begin{equation}\label{6.15}
K_7=0.
\end{equation}
Thus, for $1\le i\le m-1$ ($m\ge 2$) we have
\begin{equation}\label{6.16}
\begin{aligned}
&\lim_{t\to\infty}\frac{1}{t}\mathbb{E}\int_0^t[{\bf
I}\{Q_1(s-)=i, Q_2(s-)=j\}\\& -{\bf I}\{Q_1(s-)=i-1,
Q_2(s-)=j\}]\mbox{d}A(s)\\&=
\mu_1(i+1)\lim_{t\to\infty}\frac{1}{t}\int_0^t\mathbb{P}\{Q_1(s)=i+1,
Q_2(s)=j\}\mbox{d}s\\&
-\mu_1i\lim_{t\to\infty}\frac{1}{t}\int_0^t\mathbb{P}\{Q_1(s)=i,
Q_2(s)=j\}\mbox{d}s\\&
-\mu_2j\lim_{t\to\infty}\frac{1}{t}\int_0^t\mathbb{P}\{Q_1(s)=i,
Q_2(s)=j\}\mbox{d}s\\&
+\mu_2(j+1)\lim_{t\to\infty}\frac{1}{t}\int_0^t\mathbb{P}\{Q_1(s)=i-1,
Q_2(s)=j+1\}\mbox{d}s.
\end{aligned}
\end{equation}
(5.18) follows.

Now, consider the last case $i=m$. We have
\begin{equation}\label{6.17}
\begin{aligned}
&\mathbb{P}
\lim_{t\to\infty}\frac{1}{t}\Big(I_{m,j}(t)-I_{m,j}(0)\Big)\\
&=\mathbb{P}
\lim_{t\to\infty}\frac{1}{t}\int_0^tI_{m-1,j}(s-)\mbox{d}\widehat
A(s)\\&-\mathbb{P}
\lim_{t\to\infty}\frac{1}{t}\int_0^tI_{m,j}(s-)\mbox{d}\widehat
A(s)\\& +\mathbb{P}
\lim_{t\to\infty}\frac{1}{t}\int_0^tI_{m,j-1}(s-)\mbox{d}\widehat
A(s)\\&- \mu_1m\mathbb{P}
\lim_{t\to\infty}\frac{1}{t}\int_0^t[I_{m,j}(s)-I_{m,j-1}(s)]\mbox{d}s\\&
+\mu_2(j+1)\mathbb{P}
\lim_{t\to\infty}\frac{1}{t}\int_0^tI_{m-1,j+1}(s)\mbox{d}s\\&
+\mathbb{P}\lim_{t\to\infty}\frac{M_{m,j}(t)}{t}.
\end{aligned}
\end{equation}
The term of the left-hand side of \eqref{6.17} is equal to zero.
Therefore rewriting \eqref{6.17} as $0=K_1-K_2+K_3-K_4+K_5+K_6$
analogously to the above cases we have the following:
\begin{equation}\label{6.18}
K_1=\lim_{t\to\infty}\frac{1}{t}\mathbb{E}\int_0^t{\bf
I}\{Q_1(s-)=m-1, Q_2(s-)=j\}\mbox{d}A(s),
\end{equation}
\begin{equation}\label{6.19}
K_2=\lim_{t\to\infty}\frac{1}{t}\mathbb{E}\int_0^t{\bf
I}\{Q_1(s-)=m, Q_2(s-)=j\}\mbox{d}A(s),
\end{equation}
\begin{equation}\label{6.20}
K_3=\lim_{t\to\infty}\frac{1}{t}\mathbb{E}\int_0^t{\bf
I}\{Q_1(s-)=m, Q_2(s-)=j-1\}\mbox{d}A(s),
\end{equation}
\begin{equation}\label{6.21}
K_4=\mu_1m\lim_{t\to\infty}\frac{1}{t}\int_0^t\mathbb{P}\{Q_1(s)=m,
Q_2(s)=j\},
\end{equation}
\begin{equation}\label{6.22}
K_5=\mu_2(j+1)\lim_{t\to\infty}\frac{1}{t}\int_0^t\mathbb{
P}\{Q_1(s)=m-1, Q_2(s)=j+1\}\mbox{d}s,
\end{equation}
\begin{equation}\label{6.23}
K_6=0.
\end{equation}
Thus, for $i=m$,  we have the following:
\begin{equation}\label{6.24}
\begin{aligned}
&\lim_{t\to\infty}\frac{1}{t}\mathbb{E}\int_0^t[{\bf
I}\{Q_1(s-)=m, Q_2(s-)=j\}\\
& -{\bf I}\{Q_1(s-)=m-1, Q_2(s-)=j\}\\
&-{\bf I}\{Q_1(s-)=m, Q_2(s-)=j-1\} ]\mbox{d}A(s)\\&=
\mu_2(j+1)\lim_{t\to\infty}\frac{1}{t}\int_0^t\mathbb{P}\{Q_1(s)=m-1,
Q_2(s)=j+1\}\mbox{d}s\\&
-\mu_1m\lim_{t\to\infty}\frac{1}{t}\int_0^t\mathbb{P}\{Q_1(s)=m,
Q_2(s)=j\}.
\end{aligned}
\end{equation}
\eqref{5.19} follows. Theorem \ref{thm5.1} is proved.

\section{\bf Special cases}

Two corollaries of Theorem \ref{thm5.1} are provided below. The
first corollary is related to the special case when the process
$A(t)$ is Poisson. This case is well-known and can be found in
Chapter 2 of the book of Falin, and Templeton
[17]. The second corollary is related to the case of the `usual'
$A/M/m/\infty$ queue.

\begin{cor}\label{cor7.1}
If $A(t)$ is a Poisson processes with rate $\lambda$, then we have
the following system of equations.

(i) In the case $i=0$
\begin{equation}\label{7.1}
\lambda P_{0,j}=\mu_1 P_{1,j}-\mu_2j P_{0,j}.
\end{equation}

(ii) In the case $i=1,2,...,m-1$ ($m\ge 2$)
\begin{equation}\label{7.2}
\begin{aligned}
&\lambda(P_{i,j}-P_{i-1,j})\\ &= \mu_1(i+1)P_{i+1,j}-
\mu_1iP_{i,j}-\mu_2jP_{i,j}+\mu_2(j+1)P_{i-1,j+1}.
\end{aligned}
\end{equation}

(iii) In the case $i=m$
\begin{equation}\label{7.3}
\begin{aligned}
&\lambda(P_{m,j}-P_{m-1,j}-P_{m,j-1})\\&= \mu_2(j+1)P_{m-1,j+1}
-\mu_1mP_{m,j}
\end{aligned}
\end{equation}
Here in \eqref{7.1}-\eqref{7.3} we use the notation
$P_{i,j}=\lim_{t\to\infty}\mathbb{P}\{Q_1(t)=i$, $Q_2(t)=j\}$.
\end{cor}

\begin{proof} The proof of Corollary \ref{cor7.1} follows easily from
the statement of Theorem \ref{thm5.1}. Indeed, taking into account
\eqref{6.3} and the fact, that when the process $A(t)$ is Poisson
with rate $\lambda$, then we have $\widehat A(t)=\lambda t$.
Finally, the result follows by taking into account the existence
of the limiting stationary in time probabilities, that is
\begin{equation}\label{7.4}
\begin{aligned}
P_{i,j}&=\lim_{t\to\infty}\mathbb{P}\{Q_1(t)=i, Q_2(t)=j\}\\&
=\lim_{t\to\infty}\frac{1}{t}\int_0^t\mathbb{P}\{Q_1(s)=i,
Q_2(s)=j\}\mbox{d}s,\\& i=0,1,\ldots,m; \ j=0,1,\ldots.
\end{aligned}
\end{equation}
\end{proof}

\begin{cor}\label{cor7.2}
Let $\widetilde Q(t)$ denote the queue-length process for the
multiserver queueing system $A/M/1/\infty$. Then we have the
following system of equations.

(i) In the case $i=0$
\begin{equation}\label{7.5}
\begin{aligned}
&\lim_{t\to\infty}\frac{1}{t}\mathbb{E}\int_0^t{\bf I}\{\widetilde
Q(s-)=0 \}\mbox{d}A(s)\\& =\mu_1
\lim_{t\to\infty}\frac{1}{t}\int_0^t\mathbb{P}\{\widetilde
Q(s)=1\}\mbox{d}s.
\end{aligned}
\end{equation}

(ii) In the case $i=1,2,...,m-1$ ($m\ge 2$)
\begin{equation}\label{7.6}
\begin{aligned}
&\lim_{t\to\infty}\frac{1}{t}\mathbb{E}\int_0^t[{\bf
I}\{\widetilde Q(s-)=i\}-{\bf I}\{\widetilde
Q(s-)=i-1\}]\mbox{d}A(s)\\&=
\mu_1(i+1)\lim_{t\to\infty}\frac{1}{t}\int_0^t\mathbb{P}\{\widetilde
Q(s)=i+1\}\mbox{d}s\\&
-\mu_1i\lim_{t\to\infty}\frac{1}{t}\int_0^t\mathbb{P}\{\widetilde
Q(s)=i\}\mbox{d}s.
\end{aligned}
\end{equation}

(iii) In the case $i\ge m$
\begin{equation}\label{7.7}
\begin{aligned}
&\lim_{t\to\infty}\frac{1}{t}\mathbb{E}\int_0^t[{\bf
I}\{\widetilde Q(s-)=i \}-{\bf I}\{\widetilde
Q(s-)=i-1\}]\mbox{d}A(s)\\&
=\mu_1m\lim_{t\to\infty}\frac{1}{t}\int_0^t\big[\mathbb{P}\{\widetilde
Q(s)=i+1 \} -\mathbb{P}\{\widetilde Q(s)=i\}\big]\mbox{d}s.
\end{aligned}
\end{equation}
\end{cor}

\begin{proof} In order to prove Corollary \ref{cor7.2} notice that the
queue-length process $\widetilde Q(t)$ satisfies the following
equation
\begin{equation}\label{7.8}
\widetilde Q(t)=A(t)-\int_0^t\sum_{i=1}^\infty{\bf I}\{\widetilde
Q(s-)\ge i\}\mbox{d}\pi_i^{(1)}(s),
\end{equation}
where, as earlier, $\{\pi_i^{(1)}\}$ is a sequence of independent
Poisson processes with rate $\mu_1$. Then the proof of Corollary
\ref{cor7.2} is analogous to that of Theorem \ref{thm5.1}, and
based on the following equations:
\begin{equation}\label{7.9}
\begin{aligned}
\Delta I_i(t)&=I_{i-1}(t-)\Delta A(t)+I_{i+1}(t-)\Delta
\Pi_{i+1}^{(1)}(t)\\& +\Delta I_{i}(t-)[1-\Delta
A(t)][1-\Delta\Pi_{i}^{(1)}(t)]-I_i(t-),
\end{aligned}
\end{equation}
where $I_i(t)={\bf I}\{\widetilde Q(t)=i\}$.
\end{proof}

\section{\bf Asymptotic analysis of the system  as $\mu_2$ increases to
infinity}

In this section we study a behavior of the system as $\mu_2$
increases to infinity. Specifically we solve the following
problem.

$\bullet$ How behave
$$
\lim_{t\to\infty}\frac{1}{t}\int_0^t\mathbb{P}\{Q_1(s)=i,
Q_2(s)=j\}\mbox{d}s
$$
when $i<m$, $j\ge 1$?\\

The answer to the above question is given by the following
theorem.

\begin{thm}\label{thm8.1}
As $\mu_2\to\infty$, then for all $i=0,1,..., m-1$
\begin{equation}\label{8.1}
\sum_{j=1}^\infty j
\Big(\lim_{t\to\infty}\frac{1}{t}\int_0^t\mathbb{P}\{Q_1(s)=i,
Q_2(s)=j\}\mbox{d}s\Big)=O\big(\mu_2^{i-m}\big).
\end{equation}
\end{thm}

\begin{proof} Notice first, that by virtue of \eqref{4.17}
\begin{equation}\label{8.2}
\sum_{i=0}^m\sum_{j=1}^\infty
j\Big(\lim_{t\to\infty}\frac{1}{t}\int_0^t\mathbb{P}\{Q_1(s)=i,
Q_2(s)=j\}\mbox{d}s\Big)<\infty.
\end{equation}
Let us now start from \eqref{5.17}. Dividing this equation to
large parameter $\mu_2$ $(j\ge 1)$, we obtain
\begin{equation}\label{8.3}
\begin{aligned}
 &j\lim_{t\to\infty}\frac{1}{t}\int_0^t\mathbb{P}\{Q_1(s)=0,
Q_2(s)=j\}\mbox{d}s\\&=
\frac{\mu_1}{\mu_2}\lim_{t\to\infty}\frac{1}{t}\int_0^t\mathbb{
P}\{Q_1(s)=1, Q_2(s)=j\}\mbox{d}s\\&-
\frac{1}{\mu_2}\lim_{t\to\infty}\frac{1}{t}\mathbb{E}\int_0^t{\bf
I}\{Q_1(s-)=0, Q_2(s-)=j\}\mbox{d}A(s).
\end{aligned}
\end{equation}
Therefore, from \eqref{8.3} we obtain
\begin{equation}\label{8.4}
\begin{aligned}
&j\lim_{t\to\infty}\frac{1}{t}\int_0^t\mathbb{P}\{Q_1(s)=0,
Q_2(s)=j\}\mbox{d}s\\&\leq
\frac{C_0}{\mu_2}\lim_{t\to\infty}\frac{1}{t}\int_0^t\mathbb{
P}\{Q_1(s)=1, Q_2(s)=j\}\mbox{d}s,
\end{aligned}
\end{equation}
with an absolute constant $C_0$ satisfying $C_0\le\mu_1$.
Therefore,
\begin{equation}\label{8.5}
\begin{aligned}
&\sum_{j=1}^\infty\Big(j\lim_{t\to\infty}\frac{1}{t}\int_0^t\mathbb{P}\{Q_1(s)=0,
Q_2(s)=j\}\mbox{d}s\Big)\\&\leq\frac{C_0}{\mu_2}
\sum_{j=1}^\infty\Big(j\lim_{t\to\infty}\frac{1}{t}\int_0^t\mathbb{P}\{Q_1(s)=1,
Q_2(s)=j\}\mbox{d}s\Big).
\end{aligned}
\end{equation}
Notice now that because of \eqref{1.1}, as $\mu_2\to\infty$, the
expressions
\begin{equation}\label{8.6}
\lim_{t\to\infty}\frac{1}{t}\int_0^t\mathbb{P}\{Q_1(s)=i,
Q_2(s)=j\}\mbox{d}s,
\end{equation}
and
\begin{equation}\label{8.7}
\lim_{t\to\infty}\frac{1}{t}\mathbb{E}\int_0^t{\bf I}\{Q_1(s-)=i,
Q_2(s-)=j\}\mbox{d}A(s)
\end{equation}
are of the same order. Then considering equation \eqref{5.18}
divided as earlier to large parameter $\mu_2$ ($j\ge 1)$,
 with the aid of induction by the same
manner we obtain
\begin{equation}\label{8.8}
\begin{aligned}
&\sum_{j=1}^\infty\Big(j\lim_{t\to\infty}\frac{1}{t}\int_0^t\mathbb{P}\{Q_1(s)=i,
Q_2(s)=j\}\mbox{d}s\Big)\\&\leq
\frac{C_{i}}{\mu_2}\sum_{j=1}^\infty\Big(j\lim_{t\to\infty}\frac{1}{t}\int_0^t\mathbb{
P}\{Q_1(s)=i+1, Q_2(s)=j\}\mbox{d}s\Big),
\end{aligned}
\end{equation}
where $C_{i}$ is an absolute constant, $i=1,2,...,m-1$ ($m\ge 2$).
The statement follows.
\end{proof}

Theorem 8.1 helps us to establish the continuity theorem. Denoting
\begin{equation}\label{8.9}
J_1=\lim_{t\to\infty}\frac{1}{t}\int_0^t\mathbb{P}\{Q_1(s)=i,
Q_2(s)=j\}\mbox{d}s,
\end{equation}
and
\begin{equation}\label{8.10}
J_2=\lim_{t\to\infty}\frac{1}{t}\int_0^t\mathbb{P}\{\widetilde
Q(s)=i+j\}\mbox{d}s,
\end{equation}
we have the following.

\medskip

\begin{thm}\label{thm8.2}
Let $\mu_2\to\infty$. Then in the cases (1) $j=0$ and (2) $i=m$
the difference between $J_1$ and $J_2$ is $o(1)$.
\end{thm}

\begin{proof} It is clear that the difference between $J_1$
and $J_2$ should be studied only in the two cases mentioned in the
theorem, since in other cases, as $\mu_2\to\infty$, $J_1$ tends to
0 while $J_2$ remains positive in general.

In the case $j=0$ we have the following. When $i=0$, \eqref{5.17}
coincides with \eqref{7.5}. When $i=1,2,\ldots,m-1$, $m\ge 2$,
from \eqref{5.18} we have
\begin{equation}\label{8.11}
\begin{aligned}
&\lim_{t\to\infty}\frac{1}{t}\mathbb{E}\int_0^t[{\bf
I}\{Q_1(s-)=i, Q_2(s-)=0\}\\& -{\bf I}\{Q_1(s-)=i-1,
Q_2(s-)=0\}]\mbox{d}A(s)\\&=
\mu_1(i+1)\lim_{t\to\infty}\frac{1}{t}\int_0^t\mathbb{P}\{Q_1(s)=i+1,
Q_2(s)=0\}\mbox{d}s\\&
-\mu_1i\lim_{t\to\infty}\frac{1}{t}\int_0^t\mathbb{P}\{Q_1(s)=i,
Q_2(s)=0\}\mbox{d}s\\&
+\mu_2\lim_{t\to\infty}\frac{1}{t}\int_0^t\mathbb{P}\{Q_1(s)=i-1,
Q_2(s)=1\}\mbox{d}s.
\end{aligned}
\end{equation}
According to Theorem \ref{thm8.1} the last term of the right
hand-side vanishes as $\mu_2\to\infty$.  Therefore the limiting
relation, not containing this last term, coincides with
corresponding relation of \eqref{7.6}.

In turn, in the case $i=m$ and $j\ge 1$ from \eqref{5.19} we have
\begin{equation}\label{8.12}
\begin{aligned}
&\lim_{t\to\infty}\frac{1}{t}\mathbb{E}\int_0^t[{\bf
I}\{Q_1(s-)=m, Q_2(s-)=j\}\\& -{\bf I}\{Q_1(s-)=m,
Q_2(s-)=j-1\}]\mbox{d}A(s)\\&=
\mu_2(j+1)\lim_{t\to\infty}\frac{1}{t}\int_0^t\mathbb{P}\{Q_1(s)=m-1,
Q_2(s)=j+1\}\mbox{d}s\\&
-\mu_1m\lim_{t\to\infty}\frac{1}{t}\int_0^t\mathbb{P}\{Q_1(s)=m,
Q_2(s)=j\}\mbox{d}s+O\Big(\frac{1}{\mu_2}\Big).
\end{aligned}
\end{equation}
Denoting $Q(t)=m+Q_2(t)$, then \eqref{8.12} can be rewritten
\begin{equation}\label{8.13}
\begin{aligned}
&\lim_{t\to\infty}\frac{1}{t}\mathbb{E}\int_0^t[{\bf I}\{Q(s-)=m+j
\}-{\bf I}\{Q(s-)=m+j-1\}]\mbox{d}A(s)\\&=
\mu_2(j+1)\lim_{t\to\infty}\frac{1}{t}\int_0^t\mathbb{P}\{Q_1(s)=m-1,
Q_2(s)=j+1\}\mbox{d}s\\&
-\mu_1m\lim_{t\to\infty}\frac{1}{t}\int_0^t\mathbb{
P}\{Q(s)=m+j\}\mbox{d}s+O\Big(\frac{1}{\mu_2}\Big).
\end{aligned}
\end{equation}
Now, comparison with \eqref{7.5}-\eqref{7.7} and the normalization
condition enables us to conclude that
\begin{equation}\label{8.14}
\begin{aligned}
&\lim_{\mu_2\to\infty}
\mu_2(j+1)\lim_{t\to\infty}\frac{1}{t}\int_0^t\mathbb{P}\{Q_1(s)=m-1,
Q_2(s)=j+1\}\mbox{d}s\\&=\lim_{\mu_2\to\infty}
\lim_{t\to\infty}\frac{1}{t}\int_0^t\mathbb{
P}\{Q(s)=m+j+1\}\mbox{d}s,
\end{aligned}
\end{equation}
and $J_2-J_1=o(1)$ as $\mu_2\to\infty$. The theorem is
proved.\end{proof}

Note, that the analogue of Theorem 8.2 for the Markovian
multiserver retrial queueing system is proved in Falin, and
Templeton
[17].

\section{\bf An algorithm for numerical calculation of the model}

The aim of this section is to provide the method for calculating
\begin{equation}\label{9.1}
\lim_{t\to\infty}\frac{1}{t}\int_0^t\mathbb{P}\{Q_1(s)=i,
Q_2(s)=j\}\mbox{d}s.
\end{equation}
It is worth first noting, that approximation of \eqref{9.1} with
the aid of simulation by straightforward manner is not an
elementary problem. Taking $T$ large and approximating \eqref{9.1}
by
\begin{equation}\label{9.2}
\frac{1}{T} \int_0^T\mathbb{P}\{Q_1(t)=i, Q_2(t)=j\}\mbox{d}t
\end{equation}
is not realistic. For satisfactory approximation of the integral
by sum it is necessary to take a small step $\Delta$. Then number
of terms should be very large, and the computational procedure
becomes complicated.

A more simple way is to use the Lebesgue theorem on dominated
convergence:
\begin{equation}\label{9.3}
\begin{aligned}
&\lim_{t\to\infty}\frac{1}{t}\int_0^t\mathbb{P}\{Q_1(s)=i,
Q_2(s)=j\}\mbox{d}s\\&
=\mathbb{P}\lim_{t\to\infty}\frac{1}{t}\int_0^t\mathbf{I}\{Q_1(s)=i,
Q_2(s)=j\}\mbox{d}s.
\end{aligned}
\end{equation}
In this case, taking $T$ large enough we estimate
\begin{equation}\label{9.4}
\frac{1}{T} \int_0^T\mathbf{I}\{Q_1(t)=i, Q_2(t)=j\}\mbox{d}t
\end{equation}
rather than \eqref{9.2}. The trajectories of
$\mathbf{I}\{Q_1(t)=i, Q_2(t)=j\}$ are step-wise, and therefore
the computation procedure of \eqref{9.4} with the aid of
simulation is much simpler than that of \eqref{9.2}.

On the other hand, paying attention that relations
\eqref{5.17}-\eqref{5.19} contain the terms
\begin{equation}\label{9.5}
\lim_{t\to\infty}\frac{1}{t}\mathbb{E}\int_0^t{\bf I}\{Q_1(s-)=i,
Q_2(s-)=j\}\mbox{d}A(s),
\end{equation}
then approximation of \eqref{9.5} by
\begin{equation}\label{9.6}
\frac{1}{T}\mathbb{E}\int_0^T{\bf I}\{Q_1(t-)=i,
Q_2(t-)=j\}\mbox{d}A(t)
\end{equation}
is in turn simpler than approximation of \eqref{9.4}.

Indeed, \eqref{9.5} and \eqref{9.6} are the Stieltjes-type
integrals. That is for a given realization of $A(t,\omega)$ they
can be represented as finite sum in the points of jump of
$A(t,\omega)$. Then, the symbol $\mathbb{E}$ in \eqref{9.6}
requires averaging of the results after a large number of
realizations. By the Lebesgue theorem on dominated convergence we
have
\begin{equation}\label{9.7}
\begin{aligned}
&\lim_{t\to\infty}\frac{1}{t}\mathbb{E}\int_0^t{\bf I}\{Q_1(s-)=i,
Q_2(s-)=j\}\mbox{d}A(s)\\&
=\mathbb{P}\lim_{t\to\infty}\frac{1}{t}\int_0^t{\bf I}\{Q_1(s-)=i,
Q_2(s-)=j\}\mbox{d}A(s).
\end{aligned}
\end{equation}
Therefore, for enough large $T$ it can be taken
\begin{equation}\label{9.8}
a_{i,j}(T)=\frac{1}{T}\int_0^T{\bf I}\{Q_1(t-)=i,
Q_2(t-)=j\}\mbox{d}A(t)
\end{equation}
rather than \eqref{9.6}. This means that it is sufficient only one
long-run simulating.

The main difference between the computation procedures for
\eqref{9.4} and \eqref{9.8} is the following. Whereas \eqref{9.8}
requires to compute only the number of jumps in interval $(0,T)$,
\eqref{9.4} takes also into account the lengths of the time
intervals that the process spends in states ($i$, $j$). This has
no essential significance for the algorithmic complexity of the
simulation program. However, the time intervals that the process
spends in phase states may vary in wide bounds, that is the
average time that the process spends in different states ($i$,
$j$) and ($k$, $l$) may have a large difference. As a result,
\eqref{9.4} is sensitive to these variations in the sense, that a
small error for time average in specific state can result an
essential error of \eqref{9.4}. Particularly, \eqref{9.4} is
sensitive to the behavior of the process in the boundary states
$(0, i)$, associated with the case where the main queue is empty.

Hence, from the computational point of view a necessary accuracy
for \eqref{9.8} can be achieved easier than that for \eqref{9.4}.
Thus, between two suggested approaches for approximation of
\eqref{9.1}, the approach based on simulating \eqref{9.8} with
subsequential numerical solution of the system of equations is
preferable than a more straightforward approach based merely on
simulation of \eqref{9.4}.

For this reason the computational procedure below is based on
\eqref{9.8} rather than on \eqref{9.4}.
\smallskip

As values $a_{i,j}(T)$ are calculated, \eqref{5.17}-\eqref{5.19}
can be approximated as follows.

(i) In the case $i=0$
\begin{equation}\label{9.9}
\begin{aligned}
&\mu_1 \lim_{t\to\infty}\frac{1}{t}\int_0^t\mathbb{P}\{Q_1(s)=1,
Q_2(s)=j\}\mbox{d}s\\&
-\mu_2j\lim_{t\to\infty}\frac{1}{t}\int_0^t\mathbb{P}\{Q_1(s)=0,
Q_2(s)=j\}\mbox{d}s\\&\approx a_{0,j}(T).
\end{aligned}
\end{equation}

(ii) In the case $i=1,2,\ldots,m-1$
\begin{equation}\label{9.10}
\begin{aligned}
&\mu_1(i+1)\lim_{t\to\infty}\frac{1}{t}\int_0^t\mathbb{P}\{Q_1(s)=i+1,
Q_2(s)=j\}\mbox{d}s\\&
-\mu_1i\lim_{t\to\infty}\frac{1}{t}\int_0^t\mathbb{P}\{Q_1(s)=i,
Q_2(s)=j\}\mbox{d}s\\&
-\mu_2j\lim_{t\to\infty}\frac{1}{t}\int_0^t\mathbb{P}\{Q_1(s)=i,
Q_2(s)=j\}\mbox{d}s\\&
+\mu_2(j+1)\lim_{t\to\infty}\frac{1}{t}\int_0^t\mathbb{P}\{Q_1(s)=i-1,
Q_2(s)=j+1\}\mbox{d}s\\& \approx a_{i,j}(T)-a_{i-1,j}(T).
\end{aligned}
\end{equation}

(iii) In the case $i=m$
\begin{equation}\label{9.11}
\begin{aligned}
&\mu_2(j+1)\lim_{t\to\infty}\frac{1}{t}\int_0^t\mathbb{P}\{Q_1(s)=m-1,
Q_2(s)=j+1\}\mbox{d}s\\&
-\mu_1m\lim_{t\to\infty}\frac{1}{t}\int_0^t\mathbb{P}\{Q_1(s)=m,
Q_2(s)=j\}\mbox{d}s\\& \approx
a_{m,j}(T)-a_{m-1,j}(T)-a_{m,j-1}(T),
\end{aligned}
\end{equation}
where $a_{m,-1}(T)\equiv 0$, and $a_{-1,j}(T)\equiv 0$.\\

Equations \eqref{9.9}-\eqref{9.11} are similar to those for the
Markovian system. The only difference that in the case of
Markovian system the values $\lim_{t\to\infty}t^{-1}$ $\int_0^t$
$\mathbb{P}\{Q_1(s)=i, Q_2(s)=j\}\mbox{d}s$ are replaced by
limiting stationary probabilities $P_{i,j}$ (see reference
\eqref{7.4}). Then, the traditional way to estimate \eqref{9.1} is
based on one of the known truncation methods. For example, two
different methods are described in Chapter 2 of the book of Falin,
and Templeton
[17]. Other method can be found in Artalejo, and Pozo
[6]. However, the analysis of the above non-Markovian system by
truncation method is much more difficult than in the case of
Markovian system. Reduction to truncated model implies that the
initial model should be replaced by state-dependent model. In the
case of Markovian system, the system of equations for the new
state-dependent model is based on the Chapman-Kolmogorov
equations. This new system of equations is an elementary
generalization of the initial system of equations. In the case of
non-Markovian model, reduction to truncated model leads to
cumbersome analysis, and it is not clear whether the system of
equation for truncated model is similar to its variant of the
Markovian case.

The algorithm below provides numerical results remaining in a
framework of the initial model. However, it is available only for
the systems with `well-defined' parameters, when the queue-length
in orbit is not large. For example, in the case of heavy load and
low retrial rate, a queue-length in orbit is large, and the
present method becomes unsatisfactory.

\smallskip

The algorithm contains the following two steps:\\

$\bullet$ Step 1 - initial simulation.

The first step enables us to obtain the values $a_{i,j}(T)$. These
values are then used in equations \eqref{9.9}-\eqref{9.11}. There
is the finite number of equations. For the small
$\epsilon=T^{-1}$, we define the number of equations $W$ as
$W=\max\{j: a_{i,j}(T)\geq\epsilon\}$. Then, according to
\eqref{9.8} the value $W$ is associated with the maximum index $j$
for which the functional
$$
\int_0^T{\bf I}\{Q_1(t-)=i, Q_2(t-)=j\}\mbox{d}A(t)
$$
takes a positive integer value. For $j>W$ the above functional
is equal to 0.\\

$\bullet$ Step 2 - solution of the equations.

As the values $a_{i,j}(T)$ are computed, we solve the
equations and find the desired approximations \eqref{9.2} for \eqref{9.1}.\\

Notice, that in the case where the value $W$ is large, it is
necessary to use one of truncation methods nevertheless. All
numerical results obtained in the present paper are associated
with the cases where $W$ is not large, and we do not follow the
truncation methods.

\section{\bf Numerical work}

In this section a few numerical examples for simple non-Markovian
retrial queueing systems is provided. Specifically, the examples
are provided for two different retrial queueing systems having two
servers. One of them is traditional, the $D/M/2$ retrial queueing
system. Its interarrival time is equal to 1. The load of the
system $\varrho=(2\mu_1)^{-1}$ varies, including low, medium and
high load. The set of rates in the orbit varies similarly,
including low, medium and high rates.

The other queueing system is not traditional. Interarrival times
are assumed to be correlated random variables as follows. The
first interarrival time, $\xi_1$, is uniformly distributed in
interval (0,2), and the $n+1$st interarrival time is recurrently
defined as $\xi_{n+1}=2-\xi_n$, $n\ge 1$. Thus, $\{\xi_n\}_{n\ge
1}$ is a strictly stationary and ergodic sequence of random
variables having the uniform in (0,2) distribution, and
$\mathbb{E}\xi_n=1$. The parameters $\mu_1$ and $\mu_2$ of this
system vary by the same manner as these parameters of the first
system.

The aim to consider so not standard system is the following. {\it
First,} the system with correlated and alternatively changed
interarrival times often appears in a large number of
applications, and especially in telecommunication networks. For
example, such situation can occur when there are several sources,
each of which sends messages by a constant deterministic interval.
{\it Second,} our main results are related to the case of arrival
point process with strictly stationary and ergodic increments, and
it is interesting to compare the results obtained for this not
traditional system with the corresponding results related to
standard queue with usual, say deterministic, arrival.

In turn, the numerical results, obtained for the $D/M/2$ retrial
queue with high retrial rate,  are compared with corresponding
numerical results for the 'usual' $D/M/2/\infty$ queue. The last
are obtained from the known analytic representations (e.g.
Borovkov
[9]).\\

For our convenience for two-server retrial queueing systems  we
use the following notation
\begin{equation}\label{10.1}
P_{i,j}=\lim_{t\to\infty}\frac{1}{t}\int_0^t\mathbb{P}\{Q_1(s)=i,
Q_2(s)=j\}\mbox{d}s.
\end{equation}
Note, that the limiting frequency $P_{i,j}$, given by
\eqref{10.1}, can be thought as {\it steady-state} probability. In
all our numerical experiments we take $T=100,000$, and therefore
$\epsilon=0.00001$.\\

Our numerical analysis we start from the $D/M/2$ retrial queueing
system. Table 1 is related to the case of relatively low load
($\mu_1=2.5$) and different retrial rates.

\begin{table}\label{tab1}
    \begin{center}
        \begin{tabular}{cl|ccc}\hline
        & Retrial rate & 0.1 & 1.0 & 10.0\\
($i,j$)        &    &Column 1     &Column 2     &Column 3\\
\hline (0,0) & &0.5988 & 0.6333 & 0.6343\\
(0,1) & &0.0151 & 0.0007 & -\\
(0,2) & &0.0017 & -      & -\\
(1,0) & &0.3621 & 0.3642 & 0.3644\\
(1,1) & &0.0176 & 0.0001 & -\\
(2,0) & &0.0013 & 0.0013 & 0.0013\\
(2,1) & &0.0002 & -      & -\\
 \hline
        \end{tabular}
        \caption{The values of $P_{i,j}$ for the case of relatively low load
$\mu_1$=2.5}
    \end{center}
\end{table}

In the case of low retrial rate $\mu_2=0.1$ by simulation we
obtain $W=2$. Specifically, $a_{0,2}(100,000)>0$, that is the
value $P_{0,2}$ is positive ($P_{0,2}\approx 0.0017$). Moreover,
the maximum value of column 1 is $P_{0,0}\approx 0.5988$, and the
values $P_{0,j}$ are greater than the corresponding values
$P_{i,j}$, for $i=1, 2$. This enables us to conclude that a
customer, who upon arrival goes to the orbit, continues to spend
there a long time while the main queue is empty.

In the case of medium retrial rate $\mu_2=1.0$ by simulation we
have only $W=1$, that is the orbit capacity does not increases 1
at the moment of arrival. From column 2 of Table 1 it is seen that
the values $P_{0,1}$ and $P_{1,1}$ are sufficiently small,
nevertheless $P_{0,1}>P_{1,1}$. This can be explained by effect of
low load. The most of time the server is empty, and the situation,
when a customer in orbit returns to the empty queue, is typical.

In the case of relatively high retrial rate $\mu_2=10.0$ the
simulation gives $W=0$. In column 3 of Table 1 there are only
three positive values for $P_{i,j}$ which are approximately the
same as steady state probabilities for the $D/M/2/\infty$ queueing
system.

The next table, Table 2, is associated with the case of medium
load ($\mu_1=1.0$). In this case the value of traffic parameter
$\varrho=0.5$.

\begin{table}
    \begin{center}
        \begin{tabular}{cl|ccc}\hline
        & Retrial rate & 0.1 & 1.0 & 10.0\\
($i,j$)        &    &Column 1     &Column 2     &Column 3\\
\hline (0,0) & &0.0104 & 0.1258 & 0.2285\\
(0,1) & &0.0155 & 0.0971 & -\\
(0,2) & &0.0308 & 0.0348 & -\\
(0,3) & &0.0454 & 0.0155 & -\\
(0,4) & &0.0512 & 0.0018 & -\\
(1,0) & &0.0104 & 0.1632 & 0.4675\\
(1,1) & &0.0173 & 0.1471 & -\\
(1,2) & &0.0322 & 0.0838 & -\\
(1,3) & &0.0489 & 0.0097 & -\\
(1,4) & &0.0580 & 0.0041 & -\\
(2,0) & &0.0036 & 0.1883 & 0.2010\\
(2,1) & &0.0071 & 0.0642 & 0.0781\\
(2,2) & &0.0144 & 0.0276 & 0.0038\\
(2,3) & &0.0211 & 0.0158 & 0.0008\\
(2,4) & &0.0281 & 0.0024 & 0.0002\\
 \hline
        \end{tabular}
        \caption{The values of $P_{i,j}$ for the case of medium load
$\mu_1$=1.0}
    \end{center}
\end{table}


The data in column 1 of Table 2, associated with the case of
relatively low retrial rate, show that the expected queue-length
in orbit is relatively long. As the retrial rate increases, the
expected queue-length in orbit decreases. Whereas for $\mu_2=0.1$
we have $W=18$, then for $\mu_2=1.0$ we have $W=8$ and for
$\mu_2=10.0$ only $W=6$.

In the next table, Table 3, the results for retrial queue with
high retrial rate and for standard $D/M/2/\infty$ queue are
compared. For the $D/M/2/\infty$ queue it is assumed that
interarrival time is equal to 1. It is also assumed that
$\mu_1=1$. Therefore the load $\varrho=0.5$. We consider the
similar $D/M/2/\infty$ retrial queue with relatively high retrial
rate $\mu_2=10.0$, and the comparison results for these two
queueing systems are given in Table 3.

 Following the known
results for the $GI/M/m/\infty$ queue given in Borovkov
[9], Section 28, Theorem 10) and denoting
\begin{equation}\label{10.2}
 \widetilde P_i=\lim_{t\to\infty}\frac{1}{t}\int_0^t\mathbb{
P}\{\widetilde Q(s)=i\}\mbox{d}s,
\end{equation}
we have
\begin{equation}\label{10.3}
\widetilde P_i =\frac{\lambda\widetilde p_{i-1}}{i\mu_1},
~i=1,2,...,m-1,
\end{equation}
\begin{equation}\label{10.4}
\widetilde P_i=\frac{\lambda\widetilde p_{i-1}}{m\mu_1},
~i=m,m+1,...,
\end{equation}
where $\lambda$ is the reciprocal of the expected interarrival
time, and
\begin{equation}\label{10.5}
\widetilde p_i=\lim_{n\to\infty}\mathbb{P}\{\widetilde
Q(t_n-)=i\},
\end{equation}
$t_n$ is the moment of the $n$th jump of the point process $A(t)$
(i.e. the moment of $n$th arrival). The explicit representation
for $\widetilde p_i$ in turn can be found in Borovkov
[9], Section 28, Theorem 9 or in Bharucha-Reid
[8].

In our case we have: $\widetilde P_1=\widetilde p_0$ and
$\widetilde P_i=0.5\widetilde p_{i-1},~ i=2,3,...$, and by
normalization condition $\widetilde P_0$=0.5(1-$\widetilde p_0$).

In turn, $\widetilde p_0=U_0-U_1$, $\widetilde p_1=U_1$,
$p_i=r\varphi^{i-2}$, $i=2,3,...$; $\varphi$ is the root of
equation $\log z=2z-2$, $\varphi\approx 0.2031$,
$$
U_0=1-\frac{r}{1-\varphi},~
U_1=rC_1\Big[\frac{1}{C_2(1-\psi_2)}~\frac{2(1-\psi_2)-2}{2(1-\varphi)-2}\Big],
$$
$ \psi_1=\mbox{e}^{-1}\approx 0.3679,~\psi_2=\mbox{e}^{-2}\approx
0.1353 $, $C_1=\psi_1(1-\psi_1)^{-1}\approx 0.6110$,
        $C_2=\psi_2(1-\psi_2)^{-1}\approx 0.1565$,
$$
r=\Big[\frac{1}{1-\varphi}+\frac{2}{C_1(1-\psi_1)}~\frac{2(1-\psi_1)-1}{2(1-\varphi)-1}+
\frac{1}{C_2(1-\psi_2)}~\frac{2(1-\psi_2)-2}{2(1-\varphi)-2}\Big]^{-1}
$$
$$
\approx 0.0823.
$$
Then, $U_0\approx 0.8967$, $U_1\approx 0.3544$.

\begin{table}
    \begin{center}
        \begin{tabular}{cc||cc}\hline
        Retrial & system & Standard & system\\
        \hline
($i,j$)        & $P_{i,j}$   & $k=i+j$     & $\widetilde P_{k}$\\
\hline (0,0) & 0.2285&0 & 0.2289\\
(1,0) &0.4675 & 1 &0.5423\\
(2,0) &0.2010 & 2 &0.1772\\
(2,1) &0.0781 & 3 &0.0412\\
(2,2) &0.0038 & 4 &0.0084\\
(2,3) &0.0008 & 5 &0.0017\\
(2,4) &0.0002 & 6 &0.0003\\
 \hline
        \end{tabular}
        \caption{The values of $P_{i,j}$ and $\widetilde P_{i+j}$ for the case of medium load
and relatively high retrial rate}
    \end{center}
\end{table}

The first 7 values of  $\widetilde P_k$, $(k=0,1,...,6)$, are in
the last column of Table 3. In the second column of this table are
corresponding values of $P_{i,j}$ taken from column 3 of Table II.
The results of Table 3 are agreed with convergence Theorem
\ref{thm8.2}.






\smallskip

We study now the cases of relatively high load, $\mu_1=0.6$. The
the cases of high load lead to increasing queue-length in orbit,
and therefore the numerical analysis becomes difficult. For
example, if in addition the retrial rate is low, then the number
of equations becomes large, and only special methods of analysis
are necessary. For example, if $\mu_2=0.1$, then the initial
simulation shows that $W>50$, and the values $a_{0,0}(100,000)$,
$a_{1,0}(100,000)$ and $a_{2,0}$ are negligible. At the same time,
$a_{2,30}\approx 0.0055$, $a_{2,39}\approx 0.0331$. In the case of
medium retrial rate as $\mu_2=1.0$ the number of equations is
still large, $W=30$. Here, the values $a_{0,0}(100,000)\approx
0.0014$, $a_{1,0}(100,000)\approx 0.0049$,
$a_{2,0}(100,000)\approx 0.0167$. The maximum value
$a_{i,j}^*(100,000)$ is achieved for $i=2$ and $j=4$. Namely,
$a_{2,4}\approx 0.1265$. In the case of relatively high retrial
rate as $\mu_2=10$, by initial simulation we obtain $W=27$.
However, the values $a_{i,j}(100,000)$ decreases in $j$, and the
maximum value $a_{i,j}^*(100,000)$ is achieved for $i=2$ and
$j=0$. Namely, $a_{2,0}(100,000)\approx 0.361$. We provide Table 4
of some values when $\mu_1=0.6$ and $\mu_2=10$. The left side of
the table (columns 1 and 2)  contains the values for retrial
queue, while the right side of the table (columns 3 and 4) is
associated with the standard multiserver queue. (The results for
the standard multiserver queue in this table are obtained by the
computations analogous to that of Table 3.)

\begin{table}
    \begin{center}
        \begin{tabular}{cc||cc}\hline
        Retrial & system & Standard & system\\
        \hline
($i,j$)        & $P_{i,j}$   & $k=i+j$     & $\widetilde P_{k}$\\
\hline (0,0) & 0.0309&0 & 0.0455\\
(1,0) &0.1910 & 1 &0.2518\\
(2,0) &0.2644 & 2 &0.3101\\
(2,1) &0.1486 & 3 &0.1891\\
(2,2) &0.1021 & 4 &0.1334\\
(2,3) &0.0700 & 5 &0.0961\\
(2,4) &0.0293 & 6 &0.0311\\
 \hline
        \end{tabular}
        \caption{The values of $P_{i,j}$ and $\widetilde P_{i+j}$ for the case of relatively high load
and relatively high retrial rate}
    \end{center}
\end{table}

The numerical analysis shows that, as load becomes high, the
difference between $P_{i,j}$ and $\widetilde P_{i+j}$ increases.\\

Now we provide numerical results for the described above
non-standard retrial. Recall that the first interarrival time
$\xi_1$ is uniformly distributed in (0,2). The other interarrival
times are determined recurrently as $\xi_{n+1}=2-\xi_n$. For this
system we provide numerical example only under medium setting
$\mu_1=1$ and relatively high retrial rate $\mu_2=10$. By initial
simulation we obtain $W=4$. This value is less than in the case of
deterministic interarrival times ($W=6$). However, whereas in the
case of deterministic interarrival the values $a_{0,1}(100,000)$
and $a_{1,1}(100,000)$ were negligible, in the case of this system
$a_{0,1}(100,000)\approx 0.0683$ and $a_{1,1}(100,000)\approx
0.0382$. In Table 5 we provide the values $P_{i,j}$ computed for
the retrial system with not standard arrivals (left side of the
table) and deterministic arrivals (right side of the table).

\begin{table}
    \begin{center}
        \begin{tabular}{l||c|c}\hline
        & Not standard  & Deterministic\\
($i,j$) & arrivals      & arrivals\\
\hline
(0,0) &0.3130 &0.2285\\
(0,1) &0.0880  &-\\
(1,0) &0.3436  &0.4675\\
(1,1) &0.0420  &-\\
(2,0) &0.2027  &0.2010\\
(2,1) &0.0100  &0.0781\\
(2,2) &0.0005  &0.0038\\
(2,3) &0.0001  &0.0008\\
(2,4) & -      &0.0002\\
 \hline
        \end{tabular}
        \caption{The values of $P_{i,j}$  for the retrial systems
        with not standard and deterministic arrivals}
    \end{center}
\end{table}




The obtained results enable us to conclude, that the behavior of
system with a not standard arrival differs from that with
deterministic interarrival time. Specifically, in the case of the
system with a not standard arrival the convergence of the
abovementioned functionals to its limits seems slower than in the
case of the system with deterministic interarrival time.

\section{\bf Concluding remarks}

In this paper, an analysis of non-Markovian multiserver retrial
queueing system is provided with the aid of the theory of
martingales. The system of equations for this system is obtained,
as well as the asymptotic analysis as parameter $\mu_2$ increases
to infinity is provided. The representation for the system of
equations enables us to study the system numerically, where some
terms of the system of equation are established by simulation.
Then the system of equation is reduced to other  system of
equations, similar to that of Markovian multiserver retrial model.
The results, obtained in the paper, enable us to study not
standard models of complex telecommunication systems arising in
the real life.

\section*{\bf Acknowledgements}

The author thanks both referees for careful reading of the paper
and a number of significant comments and suggestions, helping
substantially to improve the presentation.

\end{document}